\documentclass[10pt]{article}
\usepackage{amsmath,amssymb}
\usepackage{enumerate}
\usepackage[dvipsnames]{xcolor}

\allowdisplaybreaks[1]

\newtheorem{theorem}{Theorem}[section]

\newtheorem{lemma}[theorem]{Lemma}
\newtheorem{proposition}[theorem]{Proposition}

\newtheorem{definition}[theorem]{Definition}
\newtheorem{remark}[theorem]{Remark}

\newcommand{\Z}{\mathbb{Z}}

\renewcommand{\ker}{\operatorname{Ker}}
\newcommand{\id}{\operatorname{id}}

\newcommand{\Sym}{\operatorname{Sym}}
\newcommand{\aut}{\operatorname{Aut}}

\newcommand{\ret}{\operatorname{Ret}}
\newcommand{\soc}{\operatorname{Soc}}
\newcommand{\rad}{\operatorname{rad}}

\newcommand{\Aut}{\operatorname{Aut}}

\newcommand{\Ret}{\operatorname{Ret}}
\newcommand{\gr}{\operatorname{gr}}

\newenvironment{proof}{\par\noindent{ Proof.}}{$\qed$\par\bigskip}
\newcommand{\qed}{\enspace\vrule  height6pt  width4pt  depth2pt}
\usepackage{color}

\begin{document}
\title{New simple solutions of the Yang-Baxter equation and solutions associated to simple left braces\thanks{The first author was partially
supported by the grant MINECO PID2020-113047GB-I00 (Spain). }}
\author{F. Ced\'o \and J. Okni\'{n}ski
}
\date{}

\maketitle

\vspace{30pt}
 \noindent \begin{tabular}{llllllll}
  F. Ced\'o && J. Okni\'{n}ski \\
 Departament de Matem\`atiques &&  Institute of
Mathematics \\
 Universitat Aut\`onoma de Barcelona &&   Warsaw University \\
08193 Bellaterra (Barcelona), Spain    &&  Banacha 2, 02-097 Warsaw, Poland \\
 cedo@mat.uab.cat && okninski@mimuw.edu.pl
\end{tabular}\\

\vspace{30pt}



\begin{abstract}
Involutive non-degenerate set theoretic solutions of the
Yang-Baxter equation are considered, with a focus on finite
solutions. A rich class of indecomposable and irretractable
solutions is determined and necessary and sufficient conditions
are found in order that these solutions are simple. Then a link
between simple solutions and simple left braces is established,
that allows us to construct more examples of simple solutions. In
particular, the results answer some problems stated in the recent
paper \cite{CO21}.
\end{abstract}

\section{Introduction}

The Yang-Baxter equation is one of the important equations of
mathematical physics, that originated from the papers of Yang
\cite{Yang} and Baxter \cite{Baxter}. It also plays a key role in
several purely mathematical contexts, in particular in the
foundations of quantum groups and Hopf algebras,
\cite{BrownGoodearl,K}. To find all the solutions of the
Yang-Baxter equation is a difficult and important open problem.
Drinfeld in \cite{drinfeld} suggested the study of set-theoretic
solutions of the Yang-Baxter equation. Recall that these are pairs
$(X,r)$, where $X$ is a non-empty set and $r\colon X\times X\to
X\times X$ is a map such that
$$(r\times\id_X)\circ(\id_X\times \ r)\circ (r\times\id_X)=(\id_X\times \ r)
\circ (r\times\id_X)\circ(\id_X\times \ r),$$ on $X\times X\times
X$.

Gateva-Ivanova and Van den Bergh \cite{GIVdB}, and Etingof,
Schedler and Soloviev \cite{ESS} initiated the study of an
important class of set-theoretic solutions $(X,r)$ of the
Yang-Baxter equation, the involutive and non-degenerate solutions,
introducing important algebraic structures: the structure group,
the structure monoid and the structure algebra, as well as the
associated permutation group ${\mathcal G}(X,r)$. In \cite{R07},
Rump introduced another very fruitful structures to study this
special class of solutions, called braces. In particular,
${\mathcal G}(X,r)$ carries a natural structure of a brace.

For basic background and results on involutive non-degenerate set
theoretic solutions of the Yang-Baxter equation and on braces we
refer to \cite{CedoSurvey} and its references, and \cite{BCV,
CCP,CPR,CJOabund,CJOprimit,GI18,Ballesteros,Jedl-Pilit-Zam,Rump1,RumpSurvey,rump2020,Smok,SmokSmok}.

In the study of finite (involutive non-degenerate set-theoretic)
solutions of the Yang-Baxter equation (YBE), indecomposable and
irretractable solutions play a special role, in the sense that
every finite solution  leads to indecomposable and
irretractable solutions via iterated procedures of applying the
decomposition and the retract relation, introduced by Etingof,
Schedler and Soloviev in \cite{ESS}. Both approaches are based on
the permutation group ${\mathcal G}(X,r)$. Recent results from
\cite{CCP} show that every finite indecomposable solution of the
YBE is a so called dynamical extension (introduced by Vendramin in
\cite{V}) of a simple solution. On the other hand, finite
solutions with a primitive permutation group ${\mathcal G}(X,r)$,
referred to as primitive solutions, have been shown to be simple;
in fact this class consists only of the so called permutation
solutions determined by a cyclic permutation of prime length,
\cite{CJOprimit}. Recently in \cite{CO21}, other families of
finite simple solutions have been constructed for the first time.
The aim of this paper is to continue the study of simple solutions
of the YBE. Our approach is motivated by the decomposition of $X$
into imprimitivity blocks under the action of ${\mathcal G}(X,r)$.
It also explores and applies nontrivial relations to the
class of simple left braces, recently investigated in
\cite{BCJO18,BCJO19,CJOabund}.

The paper is organized as follows. In Section \ref{prelim} we
introduce the necessary background on solutions of the YBE and on
left braces. In Section \ref{sec4}, we improve and generalize
\cite[Theorem 4.9]{CO21}. In particular, we obtain a much bigger
family of solutions of the YBE that are irretractable and
indecomposable, and also we get necessary and sufficient
conditions for such solutions to be simple. The obtained
description turns out to be easy to apply, allowing us to
construct several concrete examples. In Section \ref{alt} we give
an alternative construction of the solutions described in Section
\ref{sec4}, and we answer in the positive Question 4.5 of
\cite{CO21}, that is, we construct finite simple solutions $(X,r)$
of the YBE such that the left brace $\mathcal{G}(X,r)$ is simple.
In Section \ref{indecomp}, we focus on the link between simple
solutions and simple (left) braces. First, we prove that a finite
non-trivial simple left brace $B$ such that there exists an orbit
$X$ by the action of the lambda map  satisfying $B=\gr(X)_+$
(meaning that, as an additive group, $B$ is generated by $X$)
naturally leads to a simple solution $(X,r)$ of the YBE. In
particular, it follows that Questions 4.4 and 4.5 of \cite{CO21}
are equivalent. Furthermore, we construct for every integer $n>1$
and prime numbers $p_1,\dots ,p_n$ a simple left brace $B$ of
order $p_1^{p_2}\cdots p_{n-1}^{p_n}p_n^{p_1}$ with an orbit $X$
by the action of the lambda map such that $B=\gr(X)_+$. We also
prove that there are finite simple left braces without this
property.

\section{Preliminaries} \label{prelim}
Let $X$ be a non-empty set and  let  $r:X\times X \rightarrow
X\times X$ be a map. For $x,y\in X$ we put $r(x,y) =(\sigma_x (y),
\gamma_y (x))$. Recall that $(X,r)$ is an involutive,
non-degenerate, set-theoretic solution of the Yang-Baxter equation
if $r^2=\id$, all the maps $\sigma_x$ and $\gamma_y$ are bijective
maps from $X$ to itself and
  $$r_{12} r_{23} r_{12} =r_{23} r_{12} r_{23},$$
where $r_{12}=r\times \id_X$ and $r_{23}=\id_X\times \ r$ are maps
from $X^3$ to itself. Because $r^{2}=\id$, one easily verifies
that $\gamma_y(x)=\sigma^{-1}_{\sigma_x(y)}(x)$, for all $x,y\in
X$ (see for example \cite[Proposition~1.6]{ESS}).

\bigskip
\noindent {\bf Convention.} Throughout the paper a solution of the
YBE will mean an involutive, non-degenerate, set-theoretic
solution
of the Yang-Baxter equation. \\

A left brace is a set $B$ with two binary operations, $+$ and
$\circ$, such that $(B,+)$ is an abelian group (the additive group
of $B$), $(B,\circ)$ is a group (the multiplicative group of $B$),
and for every $a,b,c\in B$,
 \begin{eqnarray} \label{braceeq}
  a\circ (b+c)+a&=&a\circ b+a\circ c.
 \end{eqnarray}
It is easy to check that, in any left brace, the zero and the
identity coincide. In any left brace $B$ there is an action
$\lambda\colon (B,\circ)\rightarrow \aut(B,+)$, called the lambda
map of $B$, defined by $\lambda(a)=\lambda_a$ and
$\lambda_{a}(b)=-a+a\circ b$, for $a,b\in B$. We shall write
$a\circ b=ab$, for all $a,b\in B$. A trivial brace is a left brace
$B$ such that $ab=a+b$, for all $a,b\in B$, i.e. all
$\lambda_a=\id$. The socle of a left brace $B$ is
$$\soc(B)=\{ a\in B\mid ab=a+b, \mbox{ for all
}b\in B \}.$$ Note that $\soc(B)=\ker(\lambda)$, and thus it is a
normal subgroup of the multiplicative group of $B$. The solution of
the YBE associated to a left brace $B$ is $(B,r_B)$, where
$r_B(a,b)=(\lambda_a(b),\lambda_{\lambda_a(b)}^{-1}(a))$, for all
$a,b\in B$ (see \cite[Lemma~2]{CJOComm}).

A left ideal of a left brace $B$ is a subgroup $L$ of the additive
group of $B$ such that $\lambda_a(b)\in L$, for all $b\in L$ and
all $a\in B$. An ideal of a left brace $B$ is a normal subgroup
$I$ of the multiplicative group of $B$ such that $\lambda_a(b)\in
I$, for all $b\in I$ and all $a\in B$. Note that
\begin{eqnarray}\label{addmult1}
ab^{-1}&=&a-\lambda_{ab^{-1}}(b)
\end{eqnarray}
 for all $a,b\in B$, and
    \begin{eqnarray} \label{addmult2}
     &&a-b=a+\lambda_{b}(b^{-1})= a\lambda_{a^{-1}}(\lambda_b(b^{-1}))= a\lambda_{a^{-1}b}(b^{-1}),
     \end{eqnarray}
for all $a,b\in B$. Hence, every left ideal $L$ of $B$ also is a
subgroup of the multiplicative group of $B$, and every  ideal $I$
of a left brace $B$ also is a subgroup of the additive group of
$B$. For example, it is known that $\soc(B)$ is an ideal of the
left brace $B$ (see \cite[Proposition~7]{R07}).  Note that, for
every ideal $I$ of $B$, $B/I$ inherits a natural left brace
structure.

 A homomorphism of left braces is a map $f\colon
B_1\longrightarrow B_2$, where $B_1,B_2$ are left braces, such
that $f(a b)=f(a) f(b)$ and $f(a+b)=f(a)+f(b)$, for all $a,b\in
B_1$. Note that the kernel $\ker(f)$ of a homomorphism of left
braces $f\colon B_1\longrightarrow B_2$ is an ideal of $B_1$.

Recall that if $(X,r)$ is a solution of the YBE, with
$r(x,y)=(\sigma_x(y),\gamma_y(x))$, then its structure group
$G(X,r)=\gr(x\in X\mid xy=\sigma_x(y)\gamma_y(x)$ for all
$x,y\in X)$ has a natural structure of a left brace such that
$\lambda_x(y)=\sigma_x(y)$, for all $x,y\in X$. The additive group
of $G(X,r)$ is the free abelian group with basis $X$. The
permutation group $\mathcal{G}(X,r)=\gr(\sigma_x\mid x\in X)$ of
$(X,r)$ is a subgroup of the symmetric group $\Sym_X$ on $X$.  The
map $x\mapsto \sigma_x$, from $X$ to $\mathcal{G}(X,r)$ extends to a
group homomorphism $\phi: G(X,r)\longrightarrow \mathcal{G}(X,r)$
and $\ker(\phi)=\soc(G(X,r))$. Hence there is a unique structure of
a left brace on $\mathcal{G}(X,r)$ such that $\phi$ is a
homomorphism of left braces, this is the natural structure of a left
brace on $\mathcal{G}(X,r)$.

Let $(X,r)$ and $(Y,s)$ be solutions of the YBE. We write
$r(x,y)=(\sigma_x(y),\gamma_y(x))$ and
$s(t,z)=(\sigma'_t(z),\gamma'_z(t))$, for all $x,y\in X$ and $t,z\in
Y$. A homomorphism of solutions $f\colon (X,r)\longrightarrow (Y,s)$
is a map $f\colon X\longrightarrow Y$ such that
$f(\sigma_x(y))=\sigma'_{f(x)}(f(y))$ and
$f(\gamma_y(x))=\gamma'_{f(y)}(f(x))$, for all $x,y\in X$. Since
$\gamma_y(x)=\sigma^{-1}_{\sigma_x(y)}(x)$ and
$\gamma'_z(t)=(\sigma')^{-1}_{\sigma'_t(z)}(t)$, it is clear that
$f$ is a homomorphism of solutions if and only if
$f(\sigma_x(y))=\sigma'_{f(x)}(f(y))$, for all $x,y\in X$.

Note  that, by de defining relations of the structure group, every
homomorphism of solutions $f\colon (X,r)\longrightarrow (Y,s)$
extends to a unique homomorphism of groups $f\colon
G(X,r)\longrightarrow G(Y,s)$ that we also denote by $f$, and  it
is easily checked that it also is a homomorphism of left braces
and induces a homomorphism of left braces $\bar f\colon
\mathcal{G}(X,r)\longrightarrow\mathcal{G}(Y,s)$.

In \cite{ESS}, Etingof, Schedler and Soloviev introduced the
retract relation on solutions  $(X,r)$ of the YBE. This is the
binary relation $\sim$ on $X$ defined by $x\sim y$ if and only if
$\sigma_x=\sigma_y$. Then, $\sim$ is an equivalence relation and
$r$ induces a solution $\overline{r}$ on the set
$\overline{X}=X/{\sim}$. The retract of the solution $(X,r)$ is
$\Ret(X,r)=(\overline{X},\overline{r})$. Note that the natural map
$f\colon X\longrightarrow \overline{X}:x\mapsto \bar x$ is an
epimorphism of solutions from $(X,r)$ onto $\Ret(X,r)$.

Recall that a solution $(X,r)$ is  said to be irretractable if
$\sigma_x\neq \sigma_y$ for all distinct elements $x,y\in X$, that
is $(X,r)=\Ret(X,r)$; otherwise the solution $(X,r)$ is
retractable.

\section{New indecomposable and irretractable solutions and a criterion for simplicity}\label{sec4}
Let $(X,r)$ be a solution of the YBE. We say that $(X,r)$ is
indecomposable if $\mathcal{G}(X,r)$ acts transitively on $X$.

\begin{definition}
    A solution $(X,r)$ of the YBE is simple if $|X|>1$ and for every
    epimorphism $f:(X,r) \longrightarrow (Y,s)$ of solutions either $f$
    is an isomorphism or $|Y|=1$.
\end{definition}

In this section we will improve and generalize
\cite[Theorem 4.9]{CO21}. In particular, we will obtain a much bigger
family of solutions of the YBE that are irretractable and
indecomposable, and also we will get necessary and sufficient
conditions for such solutions to be simple.

Let $A$ be an abelian additive group. Let $(j_a)_{a\in A}$ be a family of elements of $A$ such that $j_a=j_{-a}$ for all $a\in A$. We define
$$r\colon A^2\times A^2\longrightarrow A^2\times A^2$$
by
$$r((a_1,a_2),(c_1,c_2))=\left(\sigma_{(a_1,a_2)}(c_1,c_2),\sigma_{\sigma_{(a_1,a_2)}(c_1,c_2)}^{-1}(a_1,a_2)\right),$$
where
$$\sigma_{(a_1,a_2)}(c_1,c_2)=(c_1+a_2, c_2-j_{c_1+a_2-a_1}),$$
for all $a_1,a_2,c_1,c_2\in A$. Note that $\sigma_{(a_1,a_2)}\in
\Sym_{A^2}$ and
$$\sigma_{(a_1,a_2)}^{-1}(c_1,c_2)=(c_1-a_2,c_2+j_{c_1-a_1}).$$
It is known from \cite[Proposition 2]{CJOComm} that
$(A^2,r)$ is a solution of the YBE if and only if
\begin{equation}\label{sol2}\sigma_{(a_1,a_2)}\sigma_{\sigma_{(a_1,a_2)}^{-1}(c_1,c_2)}=
    \sigma_{(c_1,c_2)}\sigma_{\sigma_{(c_1,c_2)}^{-1}(a_1,a_2)},\end{equation}
for all $a_1,a_2,c_1,c_2\in A$.
We have that
\begin{eqnarray*}
\lefteqn{\sigma_{(a_1,a_2)}\sigma_{\sigma_{(a_1,a_2)}^{-1}(c_1,c_2)}(x,y)}\\
&=&\sigma_{(a_1,a_2)}\sigma_{(c_1-a_2,c_2+j_{c_1-a_1})}(x,y)\\
&=&\sigma_{(a_1,a_2)}(x+c_2+j_{c_1-a_1},y-j_{x+c_2+j_{c_1-a_1}-c_1+a_2})\\
&=&(x+c_2+j_{c_1-a_1}+a_2,y-j_{x+c_2+j_{c_1-a_1}-c_1+a_2}-j_{x+c_2+j_{c_1-a_1}+a_2-a_1})
\end{eqnarray*}
and
\begin{eqnarray*}
\lefteqn{\sigma_{(c_1,c_2)}\sigma_{\sigma_{(c_1,c_2)}^{-1}(a_1,a_2)}(x,y)}\\
&=&(x+a_2+j_{a_1-c_1}+c_2,y-j_{x+a_2+j_{a_1-c_1}-a_1+c_2}-j_{x+a_2+j_{a_1-c_1}+c_2-c_1}).
\end{eqnarray*}
Since $j_a=j_{-a}$, we get that (\ref{sol2}) holds and therefore $(A^2,r)$ is a solution of the YBE.

It is clear that the sets $\{ (a,x) : x\in A \}, a\in A$, form a
system of imprimitivity blocks of the action of the group
${\mathcal G} (A^2,r)$ on $A^2$.

\begin{proposition}\label{ind2}
    The solution $(A^2,r)$ is indecomposable if and only if $\gr ( j_a\mid a\in A)=A$.
\end{proposition}
\begin{proof}
    Let $W=\gr ( j_a\mid a\in A)$. We shall prove that the orbit of
    $(0,0)$ under the action of $\mathcal{G}(A^2,r)$  on $A^2$ is $A\times W$.
    Note that $\sigma_{(a_1,a_2)}(a,w)=(a+a_2,w-j_{a+a_2-a_1})\in A\times W$,
    for all $a,a_1,a_2\in A$ and $w\in W$.
    Let $a\in A$ and $w\in W$. There exist $a_1,\dots ,a_m\in A$ and $\varepsilon_1,\dots ,\varepsilon_m\in\{ -1,1\}$ such that
    $w=\sum_{i=1}^m\varepsilon_i j_{a-a_i}$. We have
    \begin{eqnarray*}
        \sigma_{(a_m,0)}^{-\varepsilon_m}\cdots\sigma_{(a_1,0)}^{-\varepsilon_1}\sigma_{(0,0)}^{-1}\sigma_{(0,a)}(0,0)&=&
        \sigma_{(a_m,0)}^{-\varepsilon_m}\cdots\sigma_{(a_1,0)}^{-\varepsilon_1}\sigma_{(0,0)}^{-1}(a,-j_a)\\
        &=& \sigma_{(a_m,0)}^{-\varepsilon_m}\cdots\sigma_{(a_1,0)}^{-\varepsilon_1}(a,0)\\
        &=& \sigma_{(a_m,0)}^{-\varepsilon_m}\cdots\sigma_{(a_2,0)}^{-\varepsilon_2}(a,\varepsilon_1 j_{a-a_1})=\dots\\
        &=& \sigma_{(a_m,0)}^{-\varepsilon_m}(a,\varepsilon_1 j_{a-a_1}+\dots +\varepsilon_{m-1} j_{a-a_{m-1}})\\
        &=&(a,\varepsilon_1 j_{a-a_1}+\dots +\varepsilon_{m} j_{a-a_{m}})=(a,w).
    \end{eqnarray*}
Hence the orbit of $(0,0)$ under the action of $\mathcal{G}(A^2,r)$ is $A\times W$, and the result follows.
\end{proof}

\begin{proposition}\label{irret2}
    The solution $(A^2,r)$ is irretractable if and only if for every nonzero $a\in A$ there exists $c\in A$ such that $j_c\neq j_{c+a}$.
\end{proposition}
\begin{proof}
Let $(a_1,a_2),(c_1,c_2)\in A^2$ be such that $\sigma_{(a_1,a_2)}=\sigma_{(c_1,c_2)}$. Hence
$$(x+a_2,y-j_{x+a_2-a_1})=\sigma_{(a_1,a_2)}(x,y)=\sigma_{(c_1,c_2)}(x,y)=(x+c_2,y-j_{x+c_2-c_1}),$$
for all $x,y\in A$. Thus $a_2=c_2$ and $j_{x+a_2-a_1}=j_{x+a_2-c_1}$, for all $x\in A$. Now the result follows easily.
\end{proof}

Consider the solution $(A^2,r)$ of the YBE. Let $a\in A$ be a
nonzero element. Let $V_{a,1}=\gr ( j_{c}-j_{c+a}\mid c\in
A)$. For every $i>1$,  define $V_{a,i}=V_{a,i-1}+\gr (
j_c-j_{c+v}\mid c\in A,\; v\in V_{a,i-1})$. Let
$V_a=\sum_{i=1}^{\infty}V_{a,i}$. Note that
$V_a=\bigcup_{i=1}^{\infty}V_{a,i}$.

The following result extends  \cite[Proposition 4.2]{CO21} to the infinite case.

\begin{lemma}\label{irret3}
    Let $(X,r)$ be a simple solution of the YBE. If $|X|$ is not prime,
    then $(X,r)$ is irretractable.
\end{lemma}
\begin{proof}
If $X$ is finite, then the result is \cite[Proposition 4.2]{CO21}.
Suppose that $X$ is infinite. By \cite[Proposition 4.1]{CO21},
$(X,r)$ is indecomposable. Suppose that $(X,r)$ is retractable.
Thus the natural epimorphism of solutions $(X,r)\longrightarrow
\Ret(X,r)$ is not an isomorphism. Since $(X,r)$ is simple, the
cardinality of $\ret(X,r)$ is $1$. Hence $\sigma_x=\sigma_y$ for
all $x,y\in X$. Now it is easy to see that $(X,r)$ is isomorphic
to the solution $(\Z, s)$, where $s(i,j)=(j+1,i-1)$, for all
$i,j\in\Z$. But this solution is not simple because for every
positive integer $n>1$ the natural projection $\Z\longrightarrow
\Z/(n)$ is an epimorphism of solutions from $(\Z, s)$ to
$(\Z/(n),\bar s)$, where $\bar s(\bar i,\bar j)=(\overline{j+1},
\overline{i-1})$. So we get a contradiction, and the result
follows.
\end{proof}
\begin{theorem}\label{necessary}
    Suppose that $A\neq \{0\}$. If the solution $(A^2,r)$ is simple,
    then $V_a=A$ for every nonzero element $a\in A$.
\end{theorem}
\begin{proof} Suppose that $(A^2,r)$ is simple. Let $a\in A$ be a nonzero element.
    Let $\pi\colon
    A\longrightarrow A/V_a$ be the natural map. Note that
    \begin{eqnarray*}\lefteqn{\sigma_{(a_1+v_1,a_2+v_2)}(c_1+v_3,c_2+v_4)}\\
        &=&(c_1+a_2+v_3+v_2, c_2+v_4-j_{c_1+a_2+v_3+v_2-a_1-v_1})\\
        &=&(c_1+a_2+v_3+v_2, c_2-j_{c_1+a_2-a_1}+v_4+j_{c_1+a_2-a_1}-j_{c_1+a_2+v_3+v_2-a_1-v_1}),\end{eqnarray*}
    for all $a_1,a_2,c_1,c_2\in A$ and $v_1,v_2,v_3,v_4\in V_a$. Since
    $ j_{c_1+a_2-a_1}-j_{c_1+a_2+v_3+v_2-a_1-v_1}\in V_a$,  the map
    $$\sigma'_{(\pi(a_1),\pi(a_2))}\colon (A/V_a)^2\longrightarrow (A/V_a)^2,$$
    defined by $\sigma'_{(\pi(a_1),\pi(a_2))}(\pi(c_1),\pi(c_2))=(\pi(c_1)+\pi(a_2),\pi(c_2)-\pi(j_{c_1+a_2-a_1}))$, is well-defined.
    Furthermore, if $r'\colon (A/V_a)^2\times (A/V_a)^2\longrightarrow (A/V_a)^2\times (A/V_a)^2$ is the map defined by
    \begin{eqnarray*}\lefteqn{r'((\pi(a_1),\pi(a_2)),(\pi(c_1),\pi(c_2)))}\\
        &=&\left(\sigma'_{(\pi(a_1),\pi(a_2))}(\pi(c_1),\pi(c_2)),\left(\sigma'_{\sigma'_{(\pi(a_1),\pi(a_2))}(\pi(c_1),\pi(c_2))}\right)^{-1}(\pi(a_1),\pi(a_2))\right),\end{eqnarray*}
    then $((A/V_a)^2,r')$ is a solution of the YBE and the map
    $f\colon A^2\longrightarrow (A/V_a)^2$, defined by
    $f(a_1,a_2)=(\pi(a_1),\pi(a_2))$, is an epimorphism of solutions.
    Since $(A^2,r)$ is simple, by Lemma \ref{irret3} $(A^2,r)$ is irretractable.
    By Proposition \ref{irret2}, $V_{a,1}\neq \{0\}$. Hence, since $(A^2,r)$ is simple, we have that $V_a=A$, and the result follows.
\end{proof}

The converse can be proved if the group $A$ is finite and
non-trivial.

\begin{theorem}\label{main}
    Suppose that $A$ is finite and non-trivial. If $V_a=A$, for every nonzero
    element $a\in A$, then $(A,r)$ is simple.
\end{theorem}

\begin{proof}
Assume that $V_a=A$ for every nonzero element $a\in A$. In
particular, $\gr ( j_a\mid a\in A)=A$, and thus $(A^2,r)$
is indecomposable by Proposition \ref{ind2}. Furthermore
$V_{a,1}=\gr ( j_{c}-j_{c+a}\mid c\in A)\neq\{0\}$ for
every nonzero element $a\in A$, and thus $(A^2,r)$ is
irretractable by Proposition \ref{irret2}. Let $f\colon
(A^2,r)\longrightarrow (Y,s)$ be an epimorphism of solutions which
is not an isomorphism. Since $(A^2,r)$ is indecomposable, we have
that $(Y,s)$ also is indecomposable. By \cite[Lemma 1]{CCP},
$|f^{-1}(y)|=|f^{-1}(y')|$ for all $y,y'\in Y$.  We write
$s(y_1,y_2)=(\sigma'_{y_1}(y_2),(\sigma'_{\sigma'_{y_1}(y_2)})^{-1}(y_1))$,
for all $y_1,y_2\in Y$.

We will use the following consequence of the definition of a
homomorphism of solutions:
$$f(\sigma_{(u,v)}(x,y)) = f(\sigma_{(u',v')}(x',y'))$$
whenever $f(u,v)=f(u',v')$  and  $f(x,y)=f(x',y')$.

Note also that
if $f(u,v)=f(u',v')$  and $f(x,y)=f(x',y')$, then
    \begin{eqnarray*}\sigma'_{f(u,v)}f(\sigma_{(u,v)}^{-1}(x,y))&=&f(\sigma_{(u,v)}(\sigma_{(u,v)}^{-1}(x,y))\\
        &=&f(x,y)= f(x',y')\\
        &=&f(\sigma_{(u',v')}(\sigma_{(u',v')}^{-1}(x',y'))\\
        &=&\sigma'_{f(u',v')}f(\sigma_{(u',v')}^{-1}(x',y')),
    \end{eqnarray*}
    and since $\sigma'_{f(u,v)}=\sigma'_{f(u',v')}$ is bijective, we get that
    \begin{eqnarray}\label{inverse} f(\sigma_{(u,v)}^{-1}(x,y)) = f(\sigma_{(u',v')}^{-1}(x',y')).
    \end{eqnarray}

Suppose first that there exist $a_1,a_2,c\in A$ such that
$f(a_1,c)=f(a_2,c)$ and $a_1\neq a_2$. Let $a=a_2-a_1$. We have
that
$$f(u+c,v-j_{u+c-a_1})=f(\sigma_{(a_1,c)}(u,v))=f(\sigma_{(a_2,c)}(u,v))=f(u+c,v-j_{u+c-a_2}),$$
for all $u,v\in A$.  Thus, it also follows that
$$f(u,v)=f(u,v+j_{u-a_1}-j_{u-a_2})=f(u,v+z(j_{u-a_1}-j_{u-a_2})),$$
for all $u,v\in A$ and $z\in\Z$.

Hence, we also
get
\begin{eqnarray*}\lefteqn{f(x+v,y-j_{x+v-u})}\\
    &=&f(\sigma_{(u,v)}(x,y))\\
    &=&f(\sigma_{(u,v+z(j_{u-a_1}-j_{u-a_2}))}(x,y-t(j_{x-a_1}-j_{x-a_2})))\\
    &=&f(x+v+z(j_{u-a_1}-j_{u-a_2}),y-t(j_{x-a_1}-j_{x-a_2})\\
    &&\qquad-j_{x+v+z(j_{u-a_1}-j_{u-a_2})-u}),
    \end{eqnarray*}
for all $u,v,x,y\in A$ and $z,t\in \Z$. In particular,
\begin{eqnarray*}f(x,y)&=&f(x+z(j_{u-a_1}-j_{u-a_2}),y+j_{x-u}-t(j_{x-v-a_1}-j_{x-v-a_2})\\
    &&\qquad-j_{x+z(j_{u-a_1}-j_{u-a_2})-u}),\end{eqnarray*}
for all $x,y,u,v\in A$ and $z,t\in \Z$. Hence, for
$z=0$ and $v=x-w$,
$$f(x,y)=f(x,y-t(j_{w-a_1}-j_{w-a_2})),$$
for all $x,y,w\in A$ and $t\in \Z$. Thus
$$f(x,y)=f(x,y+v_1),$$
For all $x,y\in A$ and $v_1\in V_{a,1}$. Suppose that $f(x,y)=f(x,y+v_n)$, for all $x,y\in A$ and all $v_n\in V_{a,n}$.
Then we have that
\begin{eqnarray*}f(u+y,v-j_{u+y-x})&=&f(\sigma_{(x,y)}(u,v))\\
    &=&f(\sigma_{(x,y+v_n)}(u,v+v_n'))\\
    &=&f(u+y+v_n,v+v_n'-j_{u+y+v_n-x}),\end{eqnarray*}
for all $x,y,u,v\in A$ and $v_n,v_n'\in V_{a,n}$.  This
implies that
$$f(u,v)=f(u+v_n,v+v_n'+j_{w}-j_{v_n+w}),$$
for all $u,v,w\in A$ and $v_n,v_n'\in V_{a,n}$. Hence we also have
\begin{eqnarray*}\lefteqn{f(x+v,y-j_{x+v-u})}\\
    &=&f(\sigma_{(u,v)}(x,y))\\
    &=&f(\sigma_{(u+v_n,v+v_n'+j_{w_1}-j_{v_n+w_1})}(x+z_n,y+z_n'+j_{w_2}-j_{z_n+w_2}))\\
    &=&f(x+z_n+v+v_n'+j_{w_1}-j_{v_n+w_1},y+z_n'+j_{w_2}-j_{z_n+w_2}\\
    &&\qquad-j_{x+z_n+v+v_n'+j_{w_1}-j_{v_n+w_1}-u-v_n}),
\end{eqnarray*}
for all $u,v,x,y\in A$ and $v_n,v_n'z_n,z_n'\in V_{a,n}$. For
$v_n=0$ and $z_n=-v_n'$, this leads to
$$f(x+v,y-j_{x+v-u})=f(x+v,y+z_n'+j_{w_2}-j_{z_n+w_2}-j_{x+v-u}),$$
for all $x,y,u,v,w_2\in A$ and $z_n,z_n'\in V_{a,n}$. Thus
$$f(x,y)=f(x,y+z_n'+j_{w_2}-j_{z_n+w_2}),$$
for all $x,y,w_2\in A$ and $z_n,z_n'\in V_{a,n}$. Hence
$$f(x,y)=f(x,y+v_{n+1}),$$
for all $x,y\in A$ and $v_{n+1}\in V_{a,n+1}$. By induction, we get that
$$f(x,y)=f(x,y+v),$$
for all $x,y\in A$ and $v\in V_{a}=A$. Therefore
$$f(x,0)=f(x,y),$$
for all $x,y\in A$. This implies that
$$f(0,0)=f(0,-j_0)=f(\sigma_{(0,0)}(0,0))=f(\sigma_{(0,x)}(0,y))=f(x,y-j_{x})=f(x,y),$$
for all $x,y\in A$. Therefore $|Y|=1$ in this case, as
desired. Hence we may assume that  $f(x,y)=f(u,v)$ and $y=v$
imply  $x=u$.  We will show that this case leads to a
contradiction.

Let $f(0,0)=y_0$. Let $f^{-1}(y_0)=\{ (a_i,c_i)\mid i=1, \dots
,m\}$, where $m>1$ and $c_i\neq c_j$ for all $i\neq j$. We
may assume that $(0,0)=(a_1,c_1)$. We have that
$$f(\sigma_{(0,0)}(0,0))=f(\sigma_{(0,0)}(a_i,c_i)),$$
for all $i=1,\dots ,m$. Since $|f^{-1}(y)|=|f^{-1}(y')|$ for all
$y,y'\in Y$, for $z=f(\sigma_{(0,0)}(0,0))$, it follows that
$$f^{-1}(z)=\{ (a_i,c_i-j_{a_i})\mid i=1,\dots ,m\}.$$
We also have, for every $u,i\in \{ 1,\dots ,m\}$,
$$f(\sigma_{(a_u,c_u)}(a_i,c_i))=f(\sigma_{(0,0)}(0,0))=z.$$
Hence, for every $u\in \{ 1,\dots ,m\}$,
$$f^{-1}(z)= \{ (a_i+c_u,c_i-j_{a_i+c_u-a_u})\mid i\in \{ 1,\dots ,m\}\},$$
and since $|f^{-1}(y)|=|f^{-1}(y')|$ for all $y,y'\in Y$,
for every $i\in \{ 1,\dots ,m\}$,
$$f^{-1}(z)= \{ (a_i+c_u,c_i-j_{a_i+c_u-a_u})\mid u\in \{ 1,\dots ,m\}\},$$
because $|\{ a_i+c_u\mid u\in \{ 1,\dots ,m\}\}|=|\{ c_u\mid u\in
\{ 1,\dots ,m\}\}|=m=|f^{-1}(z)|$. Therefore,
\begin{equation}\label{subgroup}
\{ a_i\mid i\in \{ 1,\dots ,m\}\}=\{ c_u\mid u\in \{ 1,\dots
,m\}\} \end{equation} is a subgroup of $A$. Now there exists $l\in
\{ 1,\dots ,m\}$ such that $a_l=-c_2$, and we have
$$\sigma_{(a_2,c_2)}(a_l,c_l)
=(0,c_l-j_{-a_2}).$$ Hence $f(0,c_l-j_{-a_2})=z=f(0,-j_0)$, and
therefore $c_l=j_{-a_2}-j_0$. We also have that $-a_2=a_t$ for
some $t\in\{ 1,\dots ,m\}$ and $-c_t=a_v$ for some $v\in\{ 1,
\dots , m\}$. Hence
$$\sigma_{(a_t,c_t)}(a_v,c_v)
=(0,c_v-j_{-a_t}),$$ thus $f(0,c_v-j_{-a_t})=z=f(0,-j_0)$, and
therefore $c_v=j_{-a_t}-j_0=j_{a_2}-j_0=j_{-a_2}-j_0=c_l$. Hence
$v=l$ and $c_2=-a_l=-a_v=c_t$. Thus $t=2$ and $-a_2=a_t=a_2$.
Similarly one can see that $a_i=-a_i$, for all $i=2,\dots ,m$.
Suppose that $a_i=c_i$, for some $i\geq 2$. In this case,
    $$f(0,-j_0)=f(\sigma_{(0,0)}(0,0))=f(\sigma_{(a_i,c_i)}(0,0))=f(c_i,-j_{c_i-a_i})=f(c_i,-j_0),$$
and we get a contradiction, because $c_i\neq 0$. Hence $a_i\neq
c_i$, for all $i=2,\dots ,m$.  Let $H=\{ a_i\mid i\in \{ 1,\dots
,m\} \}$. In view of (\ref{subgroup}) we can write $H=\{
h_1,h_2,\dots ,h_m\}$, where $h_1=0$, $h_2=a_2$ and
$(h_2,h_3),\dots ,(h_{n-1},h_n),(h_n,h_2)\in f^{-1}(y_0)$, for $n$
such that $m\geq n\geq 3$.  Since $h_n+h_n=0$,  using
(\ref{inverse}) we get
\begin{eqnarray*}f(0,0)&=&f(\sigma_{(0,0)}^{-1}\sigma_{(0,0)}(h_{n},h_{2}))\\
    &=&f(\sigma_{(0,0)}^{-1}\sigma_{(h_{n-1},h_{n})}(h_{n},h_{2}))\\
    &=&f(\sigma_{(0,0)}^{-1}(0,h_{2}-j_{h_{n-1}}))\\
    &=&f(0,h_{2}-j_{h_{n-1}}+j_{0})\end{eqnarray*}
and
\begin{eqnarray*}f(0,0)&=&f(\sigma_{(0,0)}\sigma_{(0,0)}^{-1}(h_{n},h_{2}))\\
    &=&f(\sigma_{(h_{n-1},h_{n})}\sigma_{(0,0)}^{-1}(h_{n},h_{2}))\\
    &=&f(\sigma_{(h_{n-1},h_{n})}(h_{n},h_{2}+j_{h_{n}}))\\
    &=&f(0,h_{2}+j_{h_{n}}-j_{h_{n-1}}).\end{eqnarray*}
 Hence
$h_{2}=j_{h_{n-1}}-j_0=j_{h_{n-1}}-j_{h_{n}}$, and therefore
$j_{h_n}=j_0$. Applying a similar argument, we also have
\begin{eqnarray*}f(0,0)&=&f(\sigma_{(0,0)}^{-1}\sigma_{(h_{n},h_{2})}(h_{2},h_{3}))\\
    &=&f(\sigma_{(0,0)}^{-1}(0,h_{3}-j_{h_{n}}))
=f(0,h_{3}-j_{h_{n}}+j_{0})\end{eqnarray*}
and
\begin{eqnarray*}f(0,0)&=&f(\sigma_{(h_{n},h_{2})}\sigma_{(0,0)}^{-1}(h_{2},h_{3}))\\
    &=&f(\sigma_{(h_{n},h_{2})}(h_{2},h_{3}+j_{h_{2}}))
=f(0,h_{3}+j_{h_{2}}-j_{h_{n}}).\end{eqnarray*} Hence
$h_{3}=j_{h_{n}}-j_0=j_{h_{n}}-j_{h_{2}}$, and therefore
$j_{h_2}=j_0$. Furthermore,  if $n\geq 4$ then
\begin{eqnarray*}f(0,0)&=&f(\sigma_{(0,0)}^{-1}\sigma_{(h_i,h_{i+1})}(h_{i+1},h_{i+2}))\\
    &=&f(\sigma_{(0,0)}^{-1}(0,h_{i+2}-j_{h_i}))=f(0,h_{i+2}-j_{h_i}+j_{0})\end{eqnarray*}
and
\begin{eqnarray*}f(0,0)&=&f(\sigma_{(h_i,h_{i+1})}\sigma_{(0,0)}^{-1}(h_{i+1},h_{i+2}))\\
    &=&f(\sigma_{(h_i,h_{i+1})}(h_{i+1},h_{i+2}+j_{h_{i+1}}))
=f(0,h_{i+2}+j_{h_{i+1}}-j_{h_i}).\end{eqnarray*} Hence
$h_{i+2}=j_{h_i}-j_0=j_{h_i}-j_{h_{i+1}}$, and thus
$j_{h_{i+1}}=j_0$, for all $2\leq i\leq n-2$. Therefore
    $$h_k=h_0,$$
for all $2\leq k\leq n$. But this implies that
$a_2=h_2=j_{h_{n-1}}-j_0=0$, a contradiction.  Therefore the
result follows.
\end{proof}

Note that Propositions \ref{ind2} and \ref{irret2}, and
Theorems \ref{necessary} and \ref{main} improve and generalize
\cite[Theorem 4.9]{CO21}.

A concrete example: for $A=\Z/(6)$, we take $j_0=0$, $j_1=j_2=j_4=j_5=2$ and $j_3=5$. In this case, we have
$$A_{1,1}=A_{5,1}=\gr ( j_{i+1}-j_i \mid i\in\Z/(6))=\gr ( 2,3)=\Z/(6),$$
$$A_{2,1}=A_{4,1}=\gr ( j_{i+2}-j_i \mid i\in\Z/(6))=\gr ( 2,3)=\Z/(6)$$
and
$$A_{3,1}=\gr ( j_{i+3}-j_i \mid i\in\Z/(6))=\gr ( 5)=\Z/(6).$$
Hence, by Theorem \ref{main}, the corresponding solution $(A^2,r)$
is simple. Note that $j_0-j_1=4$ is not invertible in $\Z/(6)$,
thus from \cite[Theorem 4.9]{CO21}, we only can deduce that this
solution is indecomposable and irretractable.

Another example: for $A=(\Z/(2))^2$, we take $j_{(0,0)}=(0,0)$,
$j_{(1,0)}=(0,1)$, $j_{(0,1)}=(1,1)$ and $j_{(1,1)}=(1,0)$. In
this case, we have
\begin{eqnarray*}A_{(1,0),1}&=&\gr ( j_{(i,j)}-j_{(i+1,j)} \mid (i,j)\in(\Z/(2))^2)=\gr ( (0,1)),\\
A_{(1,0),2}&=&\gr ( (0,1))+\gr ( j_{(i,j)}-j_{(i,j+1)} \mid (i,j)\in(\Z/(2))^2)\\
&=&\gr ( (0,1),(1,1))=(\Z/(2))^2,\\
A_{(0,1),1}&=&\gr ( j_{(i,j)}-j_{(i,j+1)} \mid (i,j)\in (\Z/(2))^2)=\gr ((1,1)),\\
A_{(0,1),2}&=&\gr ( (1,1)) + \gr ( j_{(i,j)}-j_{(i+1,j+1)} \mid (i,j)\in (\Z/(2))^2)\\
&=&\gr ( (1,1),(1,0))=(\Z/(2))^2,\\
A_{(1,1),1}&=&\gr ( j_{(i,j)}-j_{(i+1,j+1)} \mid (i,j)\in (\Z/(2))^2)=\gr ((1,0)),\\
A_{(1,1),2}&=&\gr ((1,0))+\gr ( j_{(i,j)}-j_{(i+1,j)} \mid (i,j)\in (\Z/(2))^2)\\
&=&\gr ( (1,0),(0,1))=(\Z/(2))^2.\end{eqnarray*} Hence, by Theorem
\ref{main}, the corresponding solution $(A^2,r)$ is simple. One
can check that this solution of cardinality $16$ is not isomorphic
to any simple solution corresponding to the cyclic  group
$A=\Z/(4)$. In fact, $\sigma_{((0,0),(0,0))}$ has order $2$ and
it is easy to see that the order of $\sigma_{(a,b)}$ of any simple
solution corresponding to $\Z/(4)$  must be $\geq 4$, for every
$a,b\in\Z/(4)$.

\section{An alternative construction}\label{alt}
In this section we give an alternative construction of the
solutions $(A^2,r)$, described in Section \ref{sec4}, in the case
where $A$ is the additive group of a commutative ring. Note that
every finite abelian additive group satisfies this condition. We
shall see that such solutions are subsolutions $(X,r)$ of the
solution associated to a left brace $B$, in fact $X$ is invariant
by the lambda map and $\mathcal{G}(X,r)$ is isomorphic to
$\gr(X)_+/\soc (\gr(X)_+)$. We also answer Question 4.5 of
\cite{CO21} in the positive, that is, we construct finite simple
solutions $(X,r)$ of the YBE  such that the left brace
$\mathcal{G}(X,r)$ is simple.

Let $A$ be a commutative nonzero ring. Let
$$H=\{ f\colon A\longrightarrow A\mid f(x)\neq 0 \mbox{ for finitely many } x\in A\}.$$
We define a sum on $H$ by $(f+g)(x)=f(x)+g(x)$, for all $f,g\in H$
and all $x\in A$. Let $\alpha\colon (A,+)\longrightarrow
\Aut(H,+)$ be the homomorphism of groups defined by
$\alpha(a)=\alpha_a$ and $\alpha_a(f)(x)=f(x-a)$, for all $f\in H$
and all $a,x\in A$. Consider the semidirect product
$$(H,+)\rtimes_{\alpha}(A,+).$$
Let $(j_a)_{a\in A}$ be a family of elements of $A$ such that
$j_a=j_{-a}$ for all $a\in A$. We define
$$b_{(j_a)_{a\in A}}\colon H\times H\longrightarrow A$$
by
$$b_{(j_a)_{a\in A}}(f,g)=\sum_{x,y\in A}f(x)j_{x-y}g(y),$$
for all $f,g\in H$. Note that $b_{(j_a)_{a\in A}}$ is well-defined
because $f(x)\neq 0$ and $g(y)\neq 0$, for finitely many $x,y\in
A$. Note that
\begin{eqnarray*}b_{(j_a)_{a\in A}}(\alpha_c(f),\alpha_c(g))&=&\sum_{x,y\in A}\alpha_c(f)(x)j_{x-y}\alpha_c(g)(y)\\
    &=&\sum_{x,y\in A}f(x-c)j_{x-y}g(y-c)\\
    &=&\sum_{x,y\in A}f(x-c)j_{x-c-(y-c)}g(y-c)\\
    &=&\sum_{x,y\in A}f(x)j_{x-y}g(y)=b_{(j_a)_{a\in A}}(f,g),
\end{eqnarray*}
for all $f,g\in H$ and all $c\in A$. Therefore we can consider the
asymmetric product (\cite{CCS}) of the trivial braces $H$ and $A$
defined via $\alpha$ and $b_{(j_a)_{a\in A}}$:
$$B_{(j_a)_{a\in A}}=H\rtimes_{\circ}A.$$
Hence $(B_{(j_a)_{a\in A}},\circ)=H\rtimes_{\alpha}A$ and the addition is defined by
$$(f,a)+(g,c)=(f+g,a+c+b_{(j_a)_{a\in A}}(f,g)),$$
for all $f,g\in H$ and $a,c\in A$. We know that $(B_{(j_a)_{a\in
A}},+,\circ)$ is a left brace. We denote the lambda map of
$B_{(j_a)_{a\in A}}$ by $\lambda^{(j_a)_{a\in A}}$. We have
\begin{eqnarray*} \lambda^{(j_a)_{a\in A}}_{(f,a)}(g,c)&=&(f,a)(g,c)-(f,a)\\
    &=&(f+\alpha_a(g),a+c)-(f,a)\\
    &=&(\alpha_a(g),c-b_{(j_a)_{a\in A}}(f,\alpha_a(g))).
    \end{eqnarray*}
For every $a\in A$, we define $e_a\in H$ by
$$e_a(x)=\delta_{a,x},$$
for all $x\in A$.

Then we have
$\alpha_c(e_a)(x)=e_a(x-c)=\delta_{a,x-c}=\delta_{a+c,x}=e_{a+c}(x)$.
Let $X=\{(e_a,c)\mid a,c\in A\}$. Note that
\begin{eqnarray}\label{lambda3}\lambda^{(j_a)_{a\in A}}_{(e_{a_1},a_2)}(e_{c_1},c_2)&=&(\alpha_{a_2}(e_{c_1}),c_2-b_{(j_a)_{a\in A}}(e_{a_1},\alpha_{a_2}(e_{c_1})))\notag\\
&=&(e_{c_1+a_2},c_2-\sum_{x,y\in A}e_{a_1}(x)j_{x-y}e_{c_1+a_2}(y))\notag\\
&=&(e_{c_1+a_2},c_2-j_{a_1-(c_1+a_2)}).\end{eqnarray}  Let
$$r_{(j_a)_{a\in A}}\colon X\times X\longrightarrow X\times X$$
be the map defined by
$$r_{(j_a)_{a\in A}}((e_{a_1},a_2),(e_{c_1},c_2))=(\lambda^{(j_a)_{a\in A}}_{(e_{a_1},a_2)}(e_{c_1},c_2),
(\lambda^{(j_a)_{a\in A}}_{\lambda^{(j_a)_{a\in
A}}_{(e_{a_1},a_2)}(e_{c_1},c_2)})^{-1}(e_{a_1},a_2)).$$ Note that
$r_{(j_a)_{a\in A}}$ is the restriction to $X^2$ of the solution
associated to the left brace $B_{(j_a)_{a\in A}}$ (see \cite[Lemma
2]{CJOComm}). Thus $(X,r_{(j_a)_{a\in A}})$ is a solution of the
YBE.

Consider the map $f\colon A^2\longrightarrow X$ defined by
$f(a,c)=(e_a,c)$, for all $a,c\in A$. Then by (\ref{lambda3}), $f$
is an isomorphism from the solution $(A^2,r)$, defined in Section
\ref{sec4}, to $(X,r_{(j_a)_{a\in A}})$. Recall that
$$r((a_1,a_2),(c_1,c_2))=(\sigma_{(a_1,a_2)}(c_1,c_2),\sigma_{\sigma_{(a_1,a_2)}(c_1,c_2)}^{-1}(a_1,a_2)),$$
where $\sigma_{(a_1,a_2)}(c_1,c_2)=(c_1+a_2,c_2-j_{c_1+a_2-a_1})$.

Note that $\gr ( X)_+$ is the left ideal  of the left brace
$B_{(j_a)_{a\in A}}$, in particular $\gr ( X)_+$ also is a left
brace with the operations inherited from $B_{(j_a)_{a\in A}}$. We
shall calculate the socle of $\gr ( X)_+$. Let $(f,c)\in \soc(\gr
( X)_+)$. Since
$(e_{a_1},a_2)=\lambda_{(f,c)}(e_{a_1},a_2)=(e_{a_1+c},
a_2-b_{(j_a)_{a\in A}}(f,e_{a_1+c}))$, for all $a_1,a_2\in A$, we
have that $c=0$ and $b(f,e_{a_1})=0$, for all $a_1\in A$. Hence
$f$ is in the radical of the restriction of the bilinear form
$b_{(j_a)_{a\in A}}$ to $\gr(e_a\mid a\in A)_+$.

\begin{proposition}\label{B}
    The left braces $\mathcal{G}(X,r_{(j_a)_{a\in A}})$ and $\gr ( X)_+/\soc(\gr ( X)_+)$
    are isomorphic.
\end{proposition}
\begin{proof}
    The proof is similar to the proof of \cite[Proposition 6.2]{CO21}.
\end{proof}

Recall that a non-zero left brace $B$ is simple if $\{0\}$ and $B$
are the only ideals of $B$.

\begin{remark}\label{simplepermu}
    Suppose that $A$ is finite and that the solution $(A^2,r)$, where
        $$r((a_1,a_2),(c_1,c_2))=(\sigma_{(a_1,a_2)}(c_1,c_2),\sigma_{\sigma_{(a_1,a_2)}(c_1,c_2)}^{-1}(a_1,a_2)),$$
        and $\sigma_{(a_1,a_2)}(c_1,c_2)=(c_1+a_2,c_2-j_{c_1+a_2-a_1})$, is simple. Then the
    solution $(X,r_{(j_a)_{a\in A}})$ is simple, and
 by \cite[Proposition 4.3]{CO21}, the left brace $\mathcal{G}(A^2,r)\cong\gr ( X)_+/\soc(\gr ( X)_+)$
    is simple if and only if
    $$\sigma_{(0,0)}\in\gr (\sigma_{(a_1,a_2)}-\sigma_{(c_1,c_2)}\mid a_1,a_2,c_1,c_2\in A )_+.$$
\end{remark}

Now we shall construct, for every prime $p\equiv 3 \mod 4$, a
simple solution $(X,r)$ such that $|X|=4p^2$ and
$\mathcal{G}(X,r)$ is a simple left brace, answering in the
positive \cite[Question 4.5]{CO21}.

\begin{proposition}\label{Exsimple} Let $p=2k+1$ be a prime with $k$ odd. Let $A=\Z/(2p)$ and let
    $$j_i=-j_{p-i}=-j_{p+i}=j_{2p-i}=2i-1\in A,$$
    for $i=1,\dots ,k$, and
    $$j_0=-j_p=-2\sum_{i=1}^k(-1)^{i}j_{i}.$$
    Then the solution $(A^2,r)$, where
    $$r((a_1,a_2),(c_1,c_2))=(\sigma_{(a_1,a_2)}(c_1,c_2),\sigma_{\sigma_{(a_1,a_2)}(c_1,c_2)}(a_1,a_2)),$$
    and $\sigma_{(a_1,a_2)}(c_1,c_2)=(c_1+a_2,c_2-j_{c_1+a_2-a_1})$, for all $a_1,a_2,c_1,c_2\in A$, is a simple solution and
    $\mathcal{G}(A^2,r)$ is a simple left brace.
\end{proposition}
\begin{proof} Note that
    $$j_0=-2\sum_{i=1}^k(-1)^{i}j_{i}=2j_k-2\sum_{i=1}^{k-1}(-1)^{i}(2i-1)=2(2k-1)-(k-1)2=2k.$$
For $u\in A\setminus\{0,p\}$ we have $j_u-j_0=\pm(2t-1)-2k$, for
some $t\in\{1,\dots ,k\}$. Hence, if $j_u-j_0\neq p$, then
$\gr(j_u-j_0)_+=A$ and thus $V_{u,1}=A$ in this case. Suppose that
$j_u-j_0=p$. In this case, $j_u=p+2k=-1$. This means that either
$u=p-1$ or $u=p+1$. If $u=p-1$, then $j_{u+2}-j_2=j_{p+1}-3=-4$,
and $\gr(-4,p)=\Z/(2p)=A$. Hence $V_{u,1}=A$ in this case. If
$u=p+1$, then $j_{u+1}-j_1=j_{p+2}-1=-j_2-1=-4$, and also
$V_{u,1}=A$ in this case. Note that $j_p-j_0=-4k=2$ and
$j_2-j_0=3-2k$. Therefore, $V_{p,2}=A$. By Theorem \ref{main}, the
solution $(A^2,r)$ is simple.

    From (\ref{addmult1}) and \cite[Lemma 2.1]{CJOprimit} it follows that
    \begin{eqnarray*}\sigma_{(a_1,a_2)}\sigma_{(c_1,c_2)}^{-1}&=&\sigma_{(a_1,a_2)}-\lambda_{\sigma_{(a_1,a_2)}\sigma_{(c_1,c_2)}^{-1}}(\sigma_{(c_1,c_2)})\\
        &=& \sigma_{(a_1,a_2)}-\sigma_{\sigma_{(a_1,a_2)}\sigma_{(c_1,c_2)}^{-1}(c_1,c_2)},
    \end{eqnarray*}
        and, also using the definition of the lambda map and the fact that
        $\id_{A^2}$ is the zero
        of the left brace $\mathcal{G}(A^2,r)$, we get
    \begin{eqnarray*}\sigma_{(a_1,a_2)}^{-1}\sigma_{(c_1,c_2)}&=&\sigma_{(a_1,a_2)}^{-1}+\lambda_{\sigma_{(a_1,a_2)}^{-1}}(\sigma_{(c_1,c_2)})\\
        &=&-\lambda_{\sigma_{(a_1,a_2)}^{-1}}(\sigma_{(a_1,a_2)})+\sigma_{\sigma_{(a_1,a_2)}^{-1}(c_1,c_2)}\\
        &=&-\sigma_{\sigma_{(a_1,a_2)}^{-1}(a_1,a_2)}+\sigma_{\sigma_{(a_1,a_2)}^{-1}(c_1,c_2)}.
    \end{eqnarray*}
    Therefore $\sigma_{(a_1,a_2)}\sigma_{(c_1,c_2)}^{-1},\sigma_{(a_1,a_2)}^{-1}\sigma_{(c_1,c_2)}\in\gr (\sigma_{(x_1,x_2)}-\sigma_{(y_1,y_2)}\mid x_1,x_2,y_1,y_2\in A)_+$ and the products of elements of these forms also are in $I=\gr (\sigma_{(x_1,x_2)}-\sigma_{(y_1,y_2)}\mid x_1,x_2,y_1,y_2\in A)_+$.
    In particular, since $k$ is odd,
    \begin{eqnarray}\label{inI}
    \sigma_{(-k,0)}^{(-1)^k}\cdots\sigma_{(-2,0)} \sigma_{(-1,0)}^{-1}\sigma_{(0,0)}^2\sigma_{(1,0)}^{-1}\sigma_{(2,0)}\cdots\sigma_{(k,0)}^{(-1)^k}\in I.
    \end{eqnarray}

    On the other hand, we have that
    \begin{eqnarray*}\lefteqn{\sigma_{(-k,0)}^{(-1)^k}\cdots\sigma_{(-2,0)} \sigma_{(-1,0)}^{-1}\sigma_{(0,0)}^2\sigma_{(1,0)}^{-1}\sigma_{(2,0)}\cdots\sigma_{(k,0)}^{(-1)^k}(u,v)}\\
        &=&\sigma_{(-k,0)}^{(-1)^k}\cdots\sigma_{(-2,0)} \sigma_{(-1,0)}^{-1}\sigma_{(0,0)}^2\sigma_{(1,0)}^{-1}\sigma_{(2,0)}\cdots\sigma_{(k-1,0)}^{(-1)^{k-1}}(u,v+j_{u-k})\\
        &=&\dots =\sigma_{(-k,0)}^{(-1)^k}\cdots\sigma_{(-2,0)} \sigma_{(-1,0)}^{-1}\sigma_{(0,0)}^2(u,v-\sum_{i=1}^k(-1)^{i}j_{u-i})\\
        &=&\sigma_{(-k,0)}^{(-1)^k}\cdots\sigma_{(-2,0)} \sigma_{(-1,0)}^{-1}(u,v-2j_u-\sum_{i=1}^k(-1)^{i}j_{u-i})\\
        &=&\sigma_{(-k,0)}^{(-1)^k}\cdots\sigma_{(-2,0)}(u,v+j_{u+1}-2j_u-\sum_{i=1}^k(-1)^{i}j_{u-i})\\
        &=&\cdots =(u,v-\sum_{i=1}^k(-1)^{i}j_{u+i}-2j_u-\sum_{i=1}^k(-1)^{i}j_{u-i})\\
        &=&(u,v-j_u-\sum_{i=-k}^{k}(-1)^{i}j_{u+i}).
        \end{eqnarray*}
    Note that for $u=0$,  since $j_i=j_{-i}$, the definition of $j_0$
    yields
    $$\sum_{i=-k}^{k}(-1)^{i}j_{u+i}=\sum_{i=-k}^{k}(-1)^{i}j_{i}=j_0+2\sum_{i=1}^{k}(-1)^{i}j_{i}=0.$$
    Suppose that $\sum_{i=-k}^{k}(-1)^{i}j_{u+i}=0$, for some $u\in A$.
 Then
    \begin{eqnarray*}\sum_{i=-k}^{k}(-1)^{i}j_{u+1+i}&=&(-1)^{k}j_{u+1+k}+\sum_{i=-k}^{k}(-1)^{i-1}j_{u+i}-(-1)^{-k-1}j_{u-k}\\
        &=&(-1)^{k}j_{u+1+k}-\sum_{i=-k}^{k}(-1)^{i}j_{u+i}-(-1)^{-k-1}j_{u-k}\\
        &=&(-1)^{k}j_{u+1+k}+(-1)^{k}j_{u-k}=0.
    \end{eqnarray*}
Therefore, by induction $\sum_{i=-k}^{k}(-1)^{i}j_{u+i}=0$, for
all $u\in A$. Hence
$$\sigma_{(-k,0)}^{(-1)^k}\cdots\sigma_{(-2,0)} \sigma_{(-1,0)}^{-1}\sigma_{(0,0)}^2\sigma_{(1,0)}^{-1}\sigma_{(2,0)}\cdots\sigma_{(k,0)}^{(-1)^k}(u,v)=(u,v-j_u)=\sigma_{(0,0)}(u,v),$$
for all $u,v\in A$. Thus,  (\ref{inI}) implies that
$$\sigma_{(0,0)}=\sigma_{(-k,0)}^{(-1)^k}\cdots\sigma_{(-2,0)} \sigma_{(-1,0)}^{-1}\sigma_{(0,0)}^2\sigma_{(1,0)}^{-1}\sigma_{(2,0)}\cdots\sigma_{(k,0)}^{(-1)^k}\in I.$$
By Remark~\ref{simplepermu}, $\mathcal{G}(A^2,r)$ is a simple left
brace.
\end{proof}

Note that, by Proposition \ref{B}, the simple left brace $\mathcal{G}((\Z/(2p))^2,r)$ of Proposition \ref{Exsimple} is isomorphic to
$$B_{j_0,\dots ,j_{2p-1}}/\soc(B_{j_0,\dots ,j_{2p-1}})=B_{j_0,\dots ,j_{2p-1}}/(\rad(b_{j_0,\dots ,j_{2p-1}})\times\{0\}).$$

\section{New simple solutions coming from simple braces} \label{indecomp}

In this section we show that Questions 4.4 and 4.5 of \cite{CO21}
are equivalent. Thus, both have a positive answer. We also
construct new examples of finite simple solutions $(X,r)$ of the
YBE such that $\mathcal{G}(X,r)$ is a simple left brace.

Let $(X,r)$ be a finite simple solution such that $|X|$ is not
prime. By \cite[Propositions~4.1 and~4.2]{CO21}, we know that
$(X,r)$ is indecomposable and irretractable. Consider the
permutation group $\mathcal{G}=\mathcal{G} (X,r)$. Then the map
$x\mapsto \sigma_x$ is an injective morphism of solutions from
$(X,r)$ to $(\mathcal{G},r_{\mathcal{G}})$. Furthermore, by
\cite[Lemma 2.1]{CJOprimit}, $\sigma(X)=\{\sigma_x\mid x\in X\}$ is
an orbit by the action of the lambda map in the left brace
$\mathcal{G}$ that generates the multiplicative (and the additive)
group of the left brace $\mathcal{G}$.

\begin{theorem}   \label{orbit}
    Let $B$ be a finite non-trivial simple left brace such that there
    exists an orbit $X$ by the action of the lambda map such that
    $B=\gr(X)_+$. Then the solution $(X,r)$ of the YBE, where
    $$r(x,y)=(\lambda_x(y),\lambda^{-1}_{\lambda_x(y)}(x)),$$
    for all $x,y\in X$, is a simple solution of the YBE.
\end{theorem}

\begin{proof}
    Let $f\colon (X,r)\longrightarrow (Y,s)$ be an epimorphism of
    solutions. Then $f$ induces an epimorphism of left braces $\bar
    f\colon \mathcal{G}(X,r)\longrightarrow \mathcal{G}(Y,s)$. On the
    other hand, the natural inclusion $i\colon (X,r)\longrightarrow
    (B,r_B)$ induces a homomorphism of left braces $\bar
    i\colon\mathcal{G}(X,r)\longrightarrow \mathcal{G}(B,r_B)$. By
    \cite[Proposition 7]{R07} (see also \cite[Lemma 2.2]{CO21}),
    $\mathcal{G}(B,r_B)\cong B/\soc(B)$. Since $B$ is a simple
    non-trivial left brace, $B/\soc(B)\cong B$. Thus we have a
    homomorphism of left braces
    $$\widetilde{i}\colon\mathcal{G}(X,r)\longrightarrow B$$
    such that $\widetilde{i}(\sigma_x)=x$, where $\sigma_x\in\Sym_X$ is
    the restriction of $\lambda_x$ to $X$,  for all $x\in X$.

    Let $g\colon B\longrightarrow \Sym_X$ be the homomorphism defined by
    $g(b)$ is the restriction of $\lambda_b$ to $X$, for all $b\in B$.
    Since $B=\gr(X)_+$, for every $b\in B$ there exist
    $x_1,\dots ,x_n\in X$, such that
    $$b=\sum_{i=1}^nx_i=x_1\lambda^{-1}_{x_1}(x_2)\lambda^{-1}_{x_1+x_2}(x_3)\cdots \lambda^{-1}_{x_1+\dots +x_{n-1}}(x_n).$$
    Hence
    $$\lambda_b=\lambda_{x_1}\circ \lambda_{\lambda^{-1}_{x_1}(x_2)}\circ\cdots\circ \lambda_{\lambda^{-1}_{x_1+\dots +x_{n-1}}(x_n)}.$$
    Therefore, $g(b)\in \mathcal{G}(X,r)$. Now it is easy to see that
    $\widetilde{i}\circ g=\id_B$. Thus $\bar f\circ g\colon
    B\longrightarrow \mathcal{G}(Y,s)$ is an epimorphism of left
    braces. Suppose that $f$ is not an isomorphism. There exist two
    distinct elements $x,z\in X$ such that $f(x)=f(z)$, and thus $\bar
    f(\sigma_x)=\bar f(\sigma_z)$. Since
    $\widetilde{i}(\sigma_x)=x\neq z=\widetilde{i}(\sigma_z)$, we have that $\sigma_x\neq \sigma_z$ and
    $\bar f(g(x))=\bar f(\sigma_x)=\bar f(\sigma_z)=\bar f(g(z))$.
    Since $B$ is a non-trivial simple left brace and $\ker(\bar f\circ
    g)\neq \{ 0\}$, we have that $\mathcal{G}(Y,s)=\{ 0\}$. Therefore
    $(Y,s)$ is a trivial solution. Hence $f(\lambda_x(y))=f(y)$, for
    all $x,y\in X$. Since $X$ is an orbit of $B$ by the action of the
    lambda map such that $B=\gr(X)_+$, we get that
    $f(x)=f(y)$, for all $x,y\in X$, and therefore $|Y|=1$, because
    $f$ is surjective, and the result follows.
\end{proof}

Note that this result and the comment above this result prove that
\cite[Questions 4.4 and 4.5]{CO21} are equivalent. Thus
Proposition \ref{Exsimple} answers these two questions in
affirmative. In the proof of our next result we shall construct a
bigger family of simple left braces with an orbit by the action of
the lambda map that generates the additive group of the brace.

\begin{theorem} \label{example}
    Let $n>1$ be an integer. Let $p_1, \dots ,p_n$ be distinct primes.
    There exists a non-trivial finite simple left brace $B$ of order
    $p_1^{p_2}\cdots p_{n-1}^{p_n}p_n^{p_1}$ with an orbit
    $X$ by the action of the lambda map such that $B=\gr(X)_+$. Furthermore, $|X|=p_1^2\cdots p_n^2$.
\end{theorem}

\begin{proof}
    Let $T=(\Z/(p_1))^{p_2-1}\times (\Z/(p_{n-1}))^{p_n-1}\times
    (\Z/(p_n))^{p_1-1}$  and $H=\Z/(p_1)\times \cdots \times \Z/(p_n)$
    be the trivial braces on these groups. Let $\alpha\colon
    H\longrightarrow \Aut(T,+)$ be the map defined by
    \begin{eqnarray}
        \label{alpha} \alpha(x_1,\dots,
        x_n)(\vec{u}_1,\dots,\vec{u}_n)=\left(\vec{u}_1(C_1^{x_2})^t,\dots
        ,\vec{u}_{n-1}(C_{n-1}^{x_n})^t,\vec{u}_n(C_n^{x_1})^t\right),
    \end{eqnarray}
    for all $x_i\in \Z/(p_i)$ and $\vec{u}_j\in
    (\Z/(p_j))^{p_{j+1}-1}$, for $1\leq j<n$, and $\vec{u}_n\in
    (\Z/(p_n))^{p_{1}-1}$, where
    $$C_j=\left(\begin{array}{ccccc}
        0&\ldots&\ldots&0&-1\\
        1&0&&0&-1\\
        0&\ddots&\ddots&\vdots&\vdots\\
        \vdots&\ddots&1&0&-1\\
        0&\ldots&0&1&-1\end{array}\right)\in M_{p_{j+1}-1}(\Z/(p_j))
    ,$$
    for $1\leq j<n$, and
    $$C_n=\left(\begin{array}{ccccc}
        0&\ldots&\ldots&0&-1\\
        1&0&&0&-1\\
        0&\ddots&\ddots&\vdots&\vdots\\
        \vdots&\ddots&1&0&-1\\
        0&\ldots&0&1&-1\end{array}\right)\in M_{p_{1}-1}(\Z/(p_n)).$$  Let $b\colon T\times T\longrightarrow H$
    be the map defined by
    $$b((\vec{u}_1,\dots,\vec{u}_n), (\vec{v}_1,\dots,\vec{v}_n))=\left(\vec{u}_1B_1\vec{v}_1^t,\dots,\vec{u}_nB_n\vec{v}_n^t\right),$$
    for all $\vec{u}_j,\vec{v}_j\in (\Z/(p_j))^{p_{j+1}-1}$, for
    $1\leq j<n$, and $\vec{u}_n,\vec{v}_n\in (\Z/(p_n))^{p_{1}-1}$,
    where
    $$B_j=\left(\begin{array}{cccc}
        1-p_{j+1}&1&\ldots&1\\
        1&\ddots&\ddots&\vdots\\
        \vdots&\ddots&\ddots&1\\
        1&\ldots&1&1-p_{j+1}\end{array}\right)\in M_{p_{j+1}-1}(\Z/(p_j))
    ,$$
    for $1\leq j<n$, and
    $$B_n=\left(\begin{array}{cccc}
        1-p_1&1&\ldots&1\\
        1&\ddots&\ddots&\vdots\\
        \vdots&\ddots&\ddots&1\\
        1&\ldots&1&1-p_1\end{array}\right)\in M_{p_{1}-1}(\Z/(p_n)).$$
    One can check that $C_i^tB_iC_i=B_i$, for all $1\leq i\leq n$. In
    fact,
    $$B_i=-\frac{1}{2}(I_i+C_i^tC_i+(C_i^2)^tC_i^2+\dots +(C_i^{k_i-1})^tC_i^{k_i-1}),$$
    where $I_i$ is the identity matrix and $k_j=p_{j+1}-1$, for $1\leq
    j<n$, and $k_n=p_1-1$. Since $C_i^{k_i}=I_i$, we get
    $C_i^tB_iC_i=B_i$, see \cite[pages 4899-4900]{BCJO18}. Hence
    \begin{eqnarray*}\lefteqn{b(\alpha(x_1,\dots, x_n)(\vec{u}_1,\dots,\vec{u}_n), \alpha(x_1,\dots, x_n)(\vec{v}_1,\dots,\vec{v}_n))}\\
        &=&b((\vec{u}_1,\dots,\vec{u}_n), (\vec{v}_1,\dots,\vec{v}_n)),\end{eqnarray*}
    for all $x_i\in \Z/(p_i)$ and $\vec{u}_j,\vec{v}_j\in
    (\Z/(p_j))^{p_{j+1}-1}$, for $1\leq j<n$, and
    $\vec{u}_n,\vec{v}_n\in (\Z/(p_n))^{p_{1}-1}$. Hence we can
    construct the asymmetric product $T\rtimes_{\circ}H$ via the
    action $\alpha$ and the bilinear map $b$ (see \cite[p.
    1850042-28]{BCJO19}). Furthermore
    $\det(B_j)=-(-p_{j+1})^{p_{j+1}-2}\neq 0$, for $1\leq j<n$ and
    $\det(B_n)=-(-p_{1})^{p_{1}-2}\neq 0$. Note also that $C_i$ is
    invertible and $C_i-I_i$ also is invertible. By \cite[Theorem 6.2]{BCJO19} and \cite[Theorem
    3.6]{BCJO18}, the left brace $T\rtimes_{\circ}H$ is simple. Recall
    that the multiplication in the left brace $T\rtimes_{\circ}H$ is
    defined by
    \begin{eqnarray*}\lefteqn{((\vec{u}_1,\dots,\vec{u}_n),(x_1,\dots,x_n))\circ
            ((\vec{v}_1,\dots,\vec{v}_n),(y_1,\dots,y_n))}\\
        &=&((\vec{u}_1,\dots,\vec{u}_n)+\alpha(x_1,\dots,x_n)(\vec{v}_1,\dots,\vec{v}_n),(x_1,\dots ,x_n)+(y_1,\dots ,y_n)),\end{eqnarray*}
    and the addition is defined by
    \begin{eqnarray*}\lefteqn{((\vec{u}_1,\dots,\vec{u}_n),(x_1,\dots,x_n))+
            ((\vec{v}_1,\dots,\vec{v}_n),(y_1,\dots,y_n))}\\
        &=&((\vec{u}_1,\dots,\vec{u}_n)+(\vec{v}_1,\dots,\vec{v}_n),(x_1,\dots ,x_n)+(y_1,\dots ,y_n)\\
        &&\qquad +b((\vec{u}_1,\dots,\vec{u}_n),(\vec{v}_1,\dots,\vec{v}_n))),\end{eqnarray*}
    for all $x_i\in \Z/(p_i)$ and $\vec{u}_j,\vec{v}_j\in (\Z/(p_j))^{p_{j+1}-1}$, for $1\leq j<n$, and $\vec{u}_n,\vec{v}_n\in (\Z/(p_n))^{p_{1}-1}$.
    Note that
    \begin{eqnarray}\label{lambda}
        \lefteqn{\lambda_{((\vec{u}_1,\dots,\vec{u}_n),(x_1,\dots,x_n))}((\vec{v}_1,\dots,\vec{v}_n),(y_1,\dots,y_n))}\\
        &=&(\alpha(x_1,\dots, x_n)(\vec{v}_1,\dots,\vec{v}_n),(y_1,\dots
        ,y_n)\nonumber\\
        &&\qquad -b(\alpha(x_1,\dots,
        x_n)(\vec{v}_1,\dots,\vec{v}_n),(\vec{u}_1,\dots,\vec{u}_n))).
        \nonumber
    \end{eqnarray} Let $X$ be the orbit of $((\vec{e}_{1,1},\dots
    ,\vec{e}_{n,1}),(0,\dots ,0))$ under the action of the lambda map,
    where $\vec{e}_{j,1},\dots ,\vec{e}_{j,p_{j+1}-1}$ is the standard
    basis of $(\Z/(p_j))^{p_{j+1}-1}$ for $1\leq j<n$, and
    $\vec{e}_{n,1},\dots ,\vec{e}_{n,p_{1}-1}$ is the standard basis
    of $(\Z/(p_n))^{p_{1}-1}$. Note that
    $$\vec{e}_{j,k}C_{j}^t=\vec{e}_{j,k+1} \mbox{ and } \vec{e}_{j,p_{j+1}-1}C_{j}^t=-\sum_{i=1}^{p_{j+1}-1}\vec{e}_{j,i}$$
    for $1\leq j<n$ and $1\leq k<p_{j+1}-1$, and
    $$\vec{e}_{n,k}C_{n}^t=\vec{e}_{n,k+1} \mbox{ and } \vec{e}_{n,p_{1}-1}C_{n}^t=-\sum_{i=1}^{p_{1}-1}\vec{e}_{n,i}$$
    for $1\leq k<p_{1}-1$.
    From the definition of $\alpha $ it follows
    easily that
    $$((\vec{e}_{1,k_1},\dots,\vec{e}_{n,k_n}),(0,\dots,0))\in X,$$
    for all $1\leq k_j\leq p_{j+1}-1$, for $1\leq j<n$, and for all $1\leq k_n\leq p_1-1$.

    On the other hand,
    \begin{eqnarray*}\lefteqn{\lambda_{((-\vec{e}_{1,2},0,\dots ,0),(0,\dots ,0))}((\vec{e}_{1,1},\dots ,\vec{e}_{n,1}),(0,\dots ,0))}\\
        &=&((\vec{e}_{1,1},\dots ,\vec{e}_{n,1}),(1,0,\dots ,0)),\end{eqnarray*}
    \begin{eqnarray*}\lefteqn{\lambda_{((-\vec{e}_{1,2},0,\dots ,0),(0,\dots ,0))}((\vec{e}_{1,1},\dots ,\vec{e}_{n,1}),(x_1,0,\dots ,0))}\\
        &=&((\vec{e}_{1,1},\dots ,\vec{e}_{n,1}),(x_1+1,0,\dots ,0)).\end{eqnarray*}
    Similarly one can see that
    $$((\vec{e}_{1,k_1},\dots,\vec{e}_{n,k_n}),(x_1,\dots,x_n))\in X,$$
    for all $1\leq k_j\leq p_{j+1}-1$, for $1\leq j<n$, for all $1\leq k_n\leq p_1-1$ and for all $x_i\in\Z/(p_i)$.
    Now it is easy  to see that
    $T\rtimes_{\circ}H=\gr(X)_+$. Note that in fact, one can see that the elements of $X$ are of the form
    $$((\vec{w}_1,\dots ,\vec{w}_n),(x_1,\dots, x_n)),$$
    where $\vec{w}_j\in\{\vec{e}_{j,1}(C_j^t)^{k}\mid 1\leq k\leq
    p_{j+1}\}$, for $1\leq j<n$,
    $\vec{w}_n\in\{\vec{e}_{n,1}(C_n^t)^{k}\mid 1\leq k\leq p_1\}$,
    and $x_i\in\Z/(p_i)$. Hence $|X|=p_1^2\cdots p_n^2$, and the result
    follows.
\end{proof}

Thus, by
Theorem~\ref{orbit} the above construction provides us a simple
solution $(X,r)$ of the YBE of cardinality $p_1^2\cdots p_n^2$. We
can also write $X=Y\times Z$, where $Z=\Z/(p_1)\times\cdots \times
\Z/(p_n)$ and $Y=\{ (\vec{e}_{1,1}(C_1^t)^{k_1},\dots
,\vec{e}_{n,1}(C_n^t)^{k_n})\in (\Z/(p_1))^{p_2-1}\times
(\Z/(p_{n-1})^{p_n-1}\times (\Z/(p_n))^{p_1-1}\mid 1\leq k_j\leq
p_{j+1}, \mbox{ for } 1\leq j<n, \mbox{ and } 1\leq k_n\leq p_1
\}$. Then $\lambda_{(i,j)}(k,l)=(\sigma_{j}(k),
l-b(\sigma_{j}(k),i))$ for $i,k\in Y, j,l\in Z$, where
$\sigma_{j}(k)=\alpha (x_1,\dots ,x_n)(\vec{v}_1,\dots
,\vec{v}_n)$ with $j=(x_1,\dots ,x_n), k=(\vec{v}_1,\dots
,\vec{v}_n)$. Here $\sigma_j = \alpha (j)$. So these are solutions
of the form considered in \cite{CO21}, where $d_{i,k} (l)=
l-b(k,i)$.

We will show that the simple solutions obtained here do not
satisfy the hypothesis about the parameters $j_i$ used in Theorem
4.9 in \cite{CO21}. In fact, in this result in \cite{CO21} we
assume that $j_0-j_i$ is invertible for all $i\neq 0$. We shall
see that this condition is not satisfied in the simple solutions
obtained here.

Let $m=p_1\cdots p_n$ and $\varphi\colon \Z/(m)\longrightarrow
\Z/(p_1)\times \dots\times\Z/(p_n)$ be the isomorphism of rings
defined by $\varphi([k]_m)=([k]_{p_1},\dots ,[k]_{p_n})$, for
every integer $k$, where $[k]_i$ denotes the class of $k$ modulo
$i$.

Let $\psi\colon Y\longrightarrow \Z/(m)$ be the map defined by
$$\psi(\vec{e}_{1,1}(C_1^t)^{k_1},\dots
,\vec{e}_{n,1}(C_n^t)^{k_n})=\varphi^{-1}([k_n]_{p_1},[k_1]_{p_2},\dots
,[k_{n-1}]_{p_n}).$$ Since $C_j$ has order $p_{j+1}$, for $1\leq
j<n$, and $C_n$ has order $p_1$, $\psi$ is well-defined and it is
bijective. Consider the map $f\colon Y\times Z\longrightarrow
\Z/(m)\times \Z/(m)$ defined by $f(y,z)=(\psi(y),
\varphi^{-1}(z))$, for all $(y,z)\in Y\times Z$.  From
(\ref{lambda}) and (\ref{alpha}) we know that
\begin{eqnarray}\label{lambda2}\lefteqn{\!\!\!\!\!\!\!\!\!\!\!\!\!\!\!\!\lambda_{((\vec{e}_{1,1}(C_1^t)^{k_1},\dots
            ,\vec{e}_{n,1}(C_n^t)^{k_n}),(x_1,\dots ,x_n))}((\vec{e}_{1,1}(C_1^t)^{l_1},\dots
        ,\vec{e}_{n,1}(C_n^t)^{l_n}),(y_1,\dots ,y_n))}\nonumber\\
    &=&((\vec{e}_{1,1}(C_1^t)^{l_1+x_2},\dots
    ,\vec{e}_{n-1,1}(C_{n-1}^t)^{l_{n-1}+x_n},\vec{e}_{n,1}(C_n^t)^{l_n+x_1}), \nonumber\\
    &&\quad
    (y_1-b_1(\vec{e}_{1,1}(C_1^t)^{l_1+x_2},\vec{e}_{1,1}(C_1^t)^{k_1}),
    \dots ,\nonumber\\
&&\qquad y_n-b_n(\vec{e}_{n,1}(C_n^t)^{l_n+x_1},\vec{e}_{n,1}(C_n^t)^{k_n})).
\end{eqnarray}
Note that, for $1\leq j<n$, if $[l_j+x_{j+1}]_{p_{j+1}}\neq
[k_j]_{p_{j+1}}$, then
$$b_j(\vec{e}_{j,1}(C_j^t)^{l_j+x_{j+1}},\vec{e}_{j,1}(C_j^t)^{k_j})=1,$$
and if $[l_j+x_{j+1}]_{p_{j+1}}= [k_j]_{p_{j+1}}$, then
$$b_j(\vec{e}_{j,1}(C_j^t)^{l_j+x_{j+1}},\vec{e}_{j,1}(C_j^t)^{k_j})=1-[p_{j+1}]_{p_j}.$$
Hence
$$b_j(\vec{e}_{j,1}(C_j^t)^{l_j+x_{j+1}},\vec{e}_{j,1}(C_j^t)^{k_j})=1-[p_{j+1}]_{p_j}\delta_{[l_j+x_{j+1}]_{p_{j+1}},[k_j]_{p_{j+1}}},$$
and similarly
$$b_j(\vec{e}_{n,1}(C_n^t)^{l_n+x_{1}},\vec{e}_{n,1}(C_n^t)^{k_n})=1-[p_{1}]_{p_n}\delta_{[l_n+x_{1}]_{p_{1}},[k_n]_{p_{1}}}.$$

Define for every $(i,j)\in \Z/(m)\times \Z/(m)$ the
map $\sigma_{(i,j)}\colon \Z/(m)\times \Z/(m)\longrightarrow
\Z/(m)\times \Z/(m)$ by
$$\sigma_{(i,j)}(k,l)=(k+j,l-j_{k+j-i}),$$
where
$$j_k=\varphi^{-1}([1]_{p_1}-[p_2]_{p_1}\delta_{[k]_{p_2},0},\dots ,[1]_{p_{n-1}}-[p_n]_{p_{n-1}}\delta_{[k]_{p_{n}},0},[1]_{p_n}-[p_1]_{p_n}\delta_{[k]_{p_1},0}).$$
Now, using (\ref{lambda2}) we get
\begin{eqnarray*}\lefteqn{\!\!\!\!\!\!\!\!\!\! f(\lambda_{((\vec{e}_{1,1}(C_1^t)^{k_1},\dots
            ,\vec{e}_{n,1}(C_n^t)^{k_n}),(x_1,\dots ,x_n))}((\vec{e}_{1,1}(C_1^t)^{l_1},\dots
        ,\vec{e}_{n,1}(C_n^t)^{l_n}),(y_1,\dots ,y_n)))}\\
    &=&(\varphi^{-1}([l_n]_{p_1}+x_1,[l_1]_{p_2}+x_2,\dots ,[l_{n-1}]_{p_n}+x_n),\\
    &&\quad \varphi^{-1}(y_1-(1-[p_2]_{p_1}\delta_{[l_1]_{p_2}+x_2,[k_1]_{p_2}}), \dots ,\\
    &&\qquad y_{n-1}-(1-[p_n]_{p_{n-1}}\delta_{[l_{n-1}]_{p_n}+x_n,[k_{n-1}]_{p_n}}),\\
    &&\qquad y_n-(1-[p_1]_{p_n}\delta_{[l_n]_{p_1}+x_1,[k_n]_{p_1}})))\\
    &=&
    (\varphi^{-1}(([l_n]_{p_1},[l_1]_{p_2},\dots ,[l_{n-1}]_{p_n})+(x_1,\dots ,x_n)),\\
    &&\quad \varphi^{-1}((y_1,\dots y_n)-(1-[p_2]_{p_1}\delta_{[l_1]_{p_2}+x_2,[k_1]_{p_2}}), \dots ,\\
    &&\qquad 1-[p_n]_{p_{n-1}}\delta_{[l_{n-1}]_{p_n}+x_n,[k_{n-1}]_{p_n}},
    1-[p_1]_{p_n}\delta_{[l_n]_{p_1}+x_1,[k_n]_{p_1}})),\end{eqnarray*}
and
\begin{eqnarray*}\lefteqn{\sigma_{f((\vec{e}_{1,1}(C_1^t)^{k_1},\dots
            ,\vec{e}_{n,1}(C_n^t)^{k_n}),(x_1,\dots ,x_n))}f((\vec{e}_{1,1}(C_1^t)^{l_1},\dots
        ,\vec{e}_{n,1}(C_n^t)^{l_n}),(y_1,\dots ,y_n))}\\
    &=&
    \sigma_{(\varphi^{-1}([k_n]_{p_1},[k_1]_{p_2}\dots
        ,[k_{n-1}]_{p_n}),\varphi^{-1}(x_1,\dots ,x_n))}(\varphi^{-1}([l_n]_{p_1},[l_1]_{p_2}\dots
    ,[l_{n-1}]),\\
    &&\qquad \varphi^{-1}(y_1,\dots ,y_n))\\
    &=&
    (\varphi^{-1}(([l_n]_{p_1},[l_1]_{p_2},\dots ,[l_{n-1}]_{p_n})+(x_1,\dots ,x_n)),\\
    &&\quad \varphi^{-1}((y_1,\dots y_n)-(1-[p_2]_{p_1}\delta_{[l_1]_{p_2}+x_2,[k_1]_{p_2}}), \dots ,\\
    &&\qquad 1-[p_n]_{p_{n-1}}\delta_{[l_{n-1}]_{p_n}+x_n,[k_{n-1}]_{p_n}},
    1-[p_1]_{p_n}\delta_{[l_n]_{p_1}+x_1,[k_n]_{p_1}})).\end{eqnarray*}
Hence $(\Z/(m)\times \Z/(m), r')$, where
$r'((i,j),(k,l))=(\sigma_{(i,j)}(k,l),\sigma^{-1}_{\sigma_{(i,j)}(k,l)}(i,j))$,
is  a solution of the YBE isomorphic to $(X,r)=(Y\times Z,r)$,
which is simple by Theorem~\ref{example}. Note that
$$j_0=\varphi^{-1}(1-p_2,\dots,1-p_n,1-p_1)$$
and $$j_{p_1}=\varphi^{-1}(1,\dots ,1,1-p_1).$$ Hence
$j_0-j_{p_1}\in\Z/(m)$ is not invertible.  Therefore,
simplicity of this solution does not follow from Theorem~4.9 in
\cite{CO21}.

A concrete example: for $m=6$, $p_1=2$ and $p_2=3$, we have
$j_0=\varphi^{-1}(1-3,1-2)=2$, $j_1=j_5=\varphi^{-1}(1,1)=1$,
$j_2=j_4=\varphi^{-1}(1,1-2)=5$ and $j_3=\varphi^{-1}(1-3,1)=4$.
Note that the solution
$((\Z/(6))^2, r)$
of the YBE, where
$r((i,j),(k,l))=(\sigma_{(i,j)}(k,l),\sigma^{-1}_{\sigma_{(i,j)}(k,l)}(i,j))$
and $\sigma_{(i,j)}(k,l)=(k+j,l-j_{k+j-i})$, for all
$i,j,k,l\in\Z/(6)$ is one of the simple solutions described in Proposition \ref{Exsimple}.

\begin{remark}\label{mod6}
    There exists a finite simple left brace $B$ such that for every orbit $X$
    of any element under the action of the lambda map, $\gr(X)_+\neq B$.
\end{remark}

\begin{proof}
    Let $p_1,p_2$ be distinct prime numbers. Suppose that $p_2>2$. Let
    $V_1=(\Z/(p_1))^{p_2-1}$ and $V_2=(\Z/(p_2))^{p_1-1}$. Let
    $T=V_1^2\times V_2$  and $H=\Z/(p_1)\times \Z/(p_2)$ be the
    trivial braces on these groups. Let $\alpha\colon H\longrightarrow
    \Aut(T,+)$ be the map defined by
    $$\alpha(x_1,x_2)(\vec{u}_{1,1},\vec{u}_{1,2},\vec{u}_2)=\left(\vec{u}_{1,1}(C_1^{x_2})^t, \vec{u}_{1,2}(C_1^{x_2})^t,\vec{u}_2(C_2^{x_1})^t\right),$$
    for all $x_i\in \Z/(p_i)$ and $\vec{u}_{1,1},\vec{u}_{1,2}\in V_1$ and $\vec{u}_2\in V_2$, where
    $$C_1=\left(\begin{array}{ccccc}
        0&\ldots&\ldots&0&-1\\
        1&0&&0&-1\\
        0&\ddots&\ddots&\vdots&\vdots\\
        \vdots&\ddots&1&0&-1\\
        0&\ldots&0&1&-1\end{array}\right)\in M_{p_{2}-1}(\Z/(p_1))
    $$
    and
    $$C_2=\left(\begin{array}{ccccc}
        0&\ldots&\ldots&0&-1\\
        1&0&&0&-1\\
        0&\ddots&\ddots&\vdots&\vdots\\
        \vdots&\ddots&1&0&-1\\
        0&\ldots&0&1&-1\end{array}\right)\in M_{p_{1}-1}(\Z/(p_2)).$$  Let $b\colon T\times T\longrightarrow H$
    be the map defined by
    $$b((\vec{u}_{1,1},\vec{u}_{1,2},\vec{u}_2), (\vec{v}_{1,1},\vec{v}_{1,2},\vec{v}_2))=
    \left(\vec{u}_{1,1}B_1\vec{v}_{1,1}^t+\vec{u}_{1,2}B_1\vec{v}_{1,2}^t,\vec{u}_2B_2\vec{v}_2^t\right),$$
    for all $\vec{u}_{1,1},\vec{u}_{1,2},\vec{v}_{1,1},\vec{v}_{1,2}\in V_1$ and $\vec{u}_2,\vec{v}_2\in V_2$, where
    $$B_1=\left(\begin{array}{cccc}
        1-p_{2}&1&\ldots&1\\
        1&\ddots&\ddots&\vdots\\
        \vdots&\ddots&\ddots&1\\
        1&\ldots&1&1-p_{2}\end{array}\right)\in M_{p_{2}-1}(\Z/(p_1))
    $$
    and
    $$B_2=\left(\begin{array}{cccc}
        1-p_1&1&\ldots&1\\
        1&\ddots&\ddots&\vdots\\
        \vdots&\ddots&\ddots&1\\
        1&\ldots&1&1-p_1\end{array}\right)\in M_{p_{1}-1}(\Z/(p_n)).$$
    As in the proof of Theorem \ref{example}, one can check that $C_i^tB_iC_i=B_i$, for all $i$, and
    we can construct the asymmetric product $T\rtimes_{\circ}H$ via the
    action $\alpha$ and the bilinear map $b$ (see \cite[p.
    1850042-28]{BCJO19}). Furthermore
    $\det(B_1)=-(-p_{2})^{p_{2}-2}\neq 0$ and
    $\det(B_2)=-(-p_{1})^{p_{1}-2}\neq 0$. Note also that $C_i$ is
    invertible and $C_i-I_i$ also is invertible. By \cite[Theorem 6.2]{BCJO19} and \cite[Theorem
    3.6]{BCJO18}, the left brace $T\rtimes_{\circ}H$ is simple.
    Note that
    \begin{eqnarray*}\lefteqn{\lambda_{((\vec{u}_{1,1},\vec{u}_{1,2},\vec{u}_2),(x_1,x_2))}((\vec{v}_{1,1},\vec{v}_{1,2},\vec{v}_2),(y_1,y_2))}\\
        &=&(\alpha(x_1,x_2)(\vec{v}_{1,1},\vec{v}_{1,2},\vec{v}_2),(y_1,y_2)-
        b(\alpha(x_1,x_2)(\vec{v}_{1,1},\vec{v}_{1,2},\vec{v}_2),(\vec{u}_{1,1},\vec{u}_{1,2},\vec{u}_2))).\end{eqnarray*}
    Let $X$ be the orbit of $((\vec{u}_{1,1},\vec{u}_{1,2},\vec{u}_2),(x_1,x_2))$ under the action of the lambda map,
    where $x_i\in\Z/(p_i)$, $\vec{u}_{1,1},\vec{u}_{1,2}\in V_1$ and
    $\vec{u}_{2}\in V_2$. Note that
    $$X\subseteq \{((\vec{u}_{1,1}(C_1^{y_2})^t,\vec{u}_{1,2}(C_1^{y_2})^t,\vec{u}_2(C_2^{y_1})^t),(z_1,z_2))\mid y_1,z_1\in\Z/(p_1),\; y_2,z_2\in \Z/(p_2) \}.$$
    Since
    $|\{(\vec{u}_{1,1}(C_1^{y_2})^t,\vec{u}_{1,2}(C_1^{y_2})^t)\mid
    y_2\in \Z/(p_2)\}|\leq p_2$  and
    $\dim_{\Z/(p_1)}(V_1^2)=2(p_2-1)>p_2$, we have that
    $\gr((\vec{u}_{1,1}(C_1^{y_2})^t,\vec{u}_{1,2}(C_1^{y_2})^t)\mid
    y_2\in \Z/(p_2))_+\neq V_1^2$. Therefore $\gr(X)_+\neq B$, and the result follows.
\end{proof}

\end{document}